\documentclass[a4paper,english,12pt]{article}
\usepackage[cp1251]{inputenc}
\usepackage{amssymb,amsmath,babel,amsthm,graphicx,cite,verbatim}
\usepackage{graphicx,color}
\usepackage{caption}
\newcounter{parag}
\newcommand{\sect}[1]
{\refstepcounter{parag}
\begin{center} {\bf \S\,\theparag. #1} \end{center}}

\newtheorem{theorem}{Theorem}

\newtheorem*{theorem*}{Theorem}
\newtheorem{lemma}{Lemma}[parag]

\theoremstyle{definition}
\newtheorem*{prf}{Proof}

\renewcommand{\pmod}[1]{\,(\mathrm{mod}\,#1)}

\newcommand{\Aut}{\operatorname{Aut}}
\newcommand{\Out}{\operatorname{Out}}

\newcounter{cl}

\newcommand{\myraise}[2][0.7]{\raisebox{#1\height}[0pt]{#2}}

\begin{document}

\begin{center}
\textbf{Almost recognizability by spectrum of simple exceptional groups of Lie type}

\textit{A.\,V.\,Vasil$'$ev\footnote{Sobolev Institute of Mathematics and Novosibirsk State University; vasand@math.nsc.ru; supported by RFBR, research project No.~13-01-00505}, A.\,M.\,Staroletov\footnote{Sobolev Institute of Mathematics and Novosibirsk State University; staroletov@math.nsc.ru; supported by RFBR, research project No.~12-01-31221}}

\end{center}

\begin{abstract}
The spectrum of a finite group is the set of its elements orders. Groups are said to be isospectral if their spectra coincide. For every finite simple exceptional group $L=E_7(q)$, we prove that each finite group isospectral to $L$ is isomorphic to a group $G$ squeezed between $L$ and its automorphism group, that is $L\leq G\leq \operatorname{Aut}L$; in particular, up-to isomorphism, there are only finitely many such groups. This assertion, together with a series of previously obtained results, implies that the same is true for every finite simple exceptional group except the group ${}^3D_4(2)$.
\end{abstract}

\medskip


\medskip


\bigskip

Given a finite group $G$, denote by $\omega(G)$ the {\it spectrum}
of $G$, i.\,e., the set of its element orders. Since for every element order all its divisors are also some element orders,
the spectrum is completely determined by the set $\mu(G)$ consisting of all maximal with respect to divisibility elements of $\omega(G)$.
We call groups $G$ and $H$ \textit{isospectral} if $\omega(G)=\omega(H)$.
Let $h(G)$ be the number of pairwise nonisomorphic groups isospectral to  $G$.
A group $G$ is called \textit{recognizable} (by spectrum) if $h(G)=1$, is \textit{almost recognizable} if $h(G)<\infty$, and is \textit{non-recognizable}
if $h(G)=\infty$. Since every finite group with a nontrivial normal soluble
subgroup is non-recognizable (see \cite[Corollary 4]{Shi1} and \cite[Lemma 1]{Mazurov}), of prime
interest is the recognition problem for nonabelian simple groups.
Following \cite{Kondratiev}, we call a finite nonabelian simple group
$L$ {\it quasirecognizable} if every finite group $G$ with
$\omega(G)=\omega(L)$ has a unique nonabelian composition factor
$S$ and $S\simeq L$. A finite group $G$ is called {\it recognizable among covers} if $\omega(G)\neq\omega(H)$
for any proper finite cover $H$ of $G$ ($H$ is a finite cover of $G$ if $G$ is a homomorphic image of $H$ and $H$ is finite). It is clear that if a finite nonabelian simple group $L$ is quasirecognizable and recognizable among covers simultaneously, then every finite group isospectral to $L$ is isomorphic to a group $G$ with $L\leq G\leq \operatorname{Aut}L$; in particular, $L$ is almost recognizable.

It turned out that many finite nonabelian simple groups are recognizable or at least almost recognizable. This paper concerns almost recognizability of finite simple exceptional groups of Lie
type, and our main purpose is to complete the proof of the following general assertion.

\begin{theorem}\label{t:main}
Let $L$ be a finite simple exceptional group of Lie type and $L\neq{}^3D_4(2)$. Then every finite group isospectral to $L$ is isomorphic to a finite group $G$ with $L\leq G\leq \operatorname{Aut}L$. In particular, $L$ is almost recognizable.
\end{theorem}

As shown in \cite{3D4}, the group  ${}^3D_4(2)$ is a real exception: it is non-recognizable and quasirecognizable at the same time.

In fact, Theorem~\ref{t:main} will follow from a series of known results and the quasirecognizability of groups $E_7(q)$ with $q>3$.

\begin{theorem}\label{t:E7} Let $L=E_7(q)$ where $q>3$. Then every finite group isospectral to $L$ is isomorphic to a group $G$ satisfying $L\leq G/K\leq\operatorname{Aut}L$, where $K$ is the maximal normal soluble subgroup of~$G$.
\end{theorem}

Indeed, the groups ${}^2B_2(q)$ \cite{ShSz}, ${}^2G_2(q)$ \cite{BrShRee}, ${}^2F_4(q)$ \cite{DSh}, $G_2(q)$ \cite{VasG2,VasStrG2}, $E_8(q)$ \cite{AlKonE8},
$F_4(2^m)$ \cite{MazShOth}, and $E_7(2)$, $E_7(3)$ \cite{14Kon} are proved to be recognizable. A recent result \cite{cover} shows that all finite simple exceptional groups besides ${}^3D_4(2)$  are recognizable among their covers. It follows that the quasirecognizability of groups ${}^3D_4(q)$ \cite{AlKonF4,Al}, $F_4(q)$ \cite{AlKonF4,AlKon03}, ${}^2E_6(q)$, $E_6(q)$ \cite{KonE6},
and Theorem~\ref{t:E7} yield the conclusion of Theorem~1.


\sect{Preliminaries}

Let $\pi$ be a set of primes. Given nonzero integer $n$, $\pi(n)$ stands for the set of all prime divisors of $n$ and $n_{\pi}$ denotes the $\pi$-part of $n$,
which is the largest positive divisor $d$ of $n$ with $\pi(d)\subseteq\pi$. The ratio $|n|/n_\pi$ is called the $\pi'$-part of $n$ and denoted by $n_{\pi'}$. For a finite group $G$, $\pi(G)=\pi(|G|)$ and $G$ is a $\pi$-group if $\pi(G)\subseteq\pi$.

For nonzero integers $n_1, n_2,\ldots,n_k$ we denote by $(n_1,n_2,\ldots,n_k)$ their greatest common divisor. The record $n_1~\mid~n_2$ means
that $n_1$ divides $n_2$, while $n_k~\vdots~n_{k-1}~\ldots~n_2~\vdots~n_1$ is a chain of divisibilities $n_1~\mid~n_2$, $n_2~\mid~n_3$, $\ldots$, $n_{k-1}~\mid~n_k$.

Let $a$ be an integer with $|a|>1$. If a prime $r$ is odd and coprime to $a$, then $e(r,a)$ denotes the multiplicative order of $a$ modulo $r$.
For an odd number $a$ we put $e(2,a)=1$, if $a\equiv1\pmod4$, and $e(2,a)=2$ if $a\equiv3\pmod4$. A prime $r$ is called a \emph{primitive prime divisor} of $a^i-1$ if $e(r,a)=i$.
The existence of primitive prime divisors for almost all pairs of $a$ and $i$ was established by Zsigmondy~\cite{Zsigmondy}.

\begin{lemma}[Zsigmondy]\label{l:zsigmondy}
Suppose that $a$ is an integer and $|a|>1$. Then for every
positive integer $i$, there is a prime $r$ with
$e(r,a)=i$ except for the cases, where $(a,i)\in\{(2,1),(2,6),(-2,2),(-2,3),(3,1),(-3,2)\}$.
\end{lemma}

The set of all primitive prime divisors of $a^i-1$ is denoted by $R_i(a)$, an element of this set is denoted by $r_i(a)$, moreover, if $a$ is fixed then the
notation $r_i$ is used. For $i\neq2$ the $R_i(a)$-part of $a^i-1$
is called the greatest primitive divisor of $a^i-1$ and denoted by $k_i(a)$. We set $k_2(a)=k_1(-a)$ and refer to it as the greatest primitive divisor of $a^2-1$. It is easy to verify that for fixed $a$ the numbers $k_i(a)$ are pairwise coprime for different $i$. Moreover, for odd $i$ we have $k_i(a)=k_{2i}(-a)$; in particular, $k_1(a)=k_2(-a)=|a-1|/2$
if $a\equiv3\pmod4$, and $k_1(a)=k_2(-a)=|a-1|$ otherwise. The following general formula \cite{R} expresses the greatest primitive divisor $k_i(a)$, $i>2$, in terms of $i$th cyclotomic polynomial $\Phi_i(x)$:
\begin{equation}\label{eq:ki}
k_i(a)=\frac{|\Phi_i(a)|}{(r,\Phi_{i_{\{r\}'}}(a))}\tag{$\ast$},
\end{equation}
where $r$ is the greatest prime divisor of $i$. Observe that if $i_{\{r\}'}$ does not divide $r-1$ then $(r,\Phi_{i_{\{r\}'}}(a))=1$.

It is known that $\Phi_m(x)=\prod\limits_{d|m}(x^d-1)^{\mu(m/d)}$, where $\mu(k)$ is the M\"{o}bius function.
The next well-known lemma collects helpful consequences of this formula (see, for example, \cite[Theorems~3.3.1 and 3.3.5]{Prasolov}).

\begin{lemma}\label{cyclomatic_polynomials}
$(1)$ Let $p$ be a prime. Then the following statements hold.
\vspace{5px}
\newline
$\Phi_{pm}(x)= \left\{
\begin {array}{l}
\Phi_m(x^p) \text{, if } (m,p)=p \text{;} \\
\Phi_m(x^p)/\Phi_n(x) \text{, if } (m, p)=1 \text{.}
\end{array} \right.$
\vspace{5px}
\newline$(2)$ If $m>1$ is an odd integer then $\Phi_{2m}(x)=\Phi_m(-x)$.
\newline$(3)$ $\Phi_p(x)=(x^p-1)/(x-1)$ and $\Phi_{2^k}(x)=x^{2^{k-1}}+1$.
\end{lemma}

\begin{lemma}\label{r-part} Let $a$ and $m$ be integers greater than $1$, and $\varepsilon\in\{+,-\}$.
\newline$(1)$ If an odd prime $r$ divides $\varepsilon{a}-1$, then $((\varepsilon{a})^m-1)_{\{r\}}=m_{\{r\}}(\varepsilon{a}-1)_{\{r\}}$.
\newline$(2)$ If an odd prime $r$ divides $(\varepsilon{a})^m-1$, then $r$ divides $(\varepsilon{a})^{m_{\{r\}'}}-1$.
\newline$(3)$ If $\varepsilon{a}-1$ is divisible by $4$, then $((\varepsilon{a})^m-1)_{\{2\}} = m_{\{2\}}(\varepsilon{a}-1)_{\{2\}}$.
\end{lemma}
\begin{prf}
See, for example, \cite[Chapter IX, Lemma 8.1]{Hup}.
\end{prf}

In notations of nonabelian simple groups we adhere to the following 	
agreements. Classical groups are considered as groups of Lie type and denoted accordingly. Furthermore, we use the short form $A^{\tau}_n(q)$
where $\tau\in\{+,-\}$, setting $A^+_n(q)=A_n(q)$ and $A^-_n(q)={}^2A_n(q)$. Similarly, we use the short form  $D^{\tau}_n(q)$ for orthogonal groups $D_n(q)$ and ${}^2D_n(q)$, where $\tau=+$ and $\tau=-$ respectively.
The alternating (symmetric) group of degree $n$ is denoted by $Alt_n$ ($Sym_n$ respectively).
For convenience we consider the Tits group ${}^2F_4(2)'$ together with sporadic groups which are denoted according to \cite{Atlas}.

Let $G$ be a finite group. The {\it prime graph} $GK(G)$
({\it Gruenberg~---~Kegel graph}) of $G$ is defined as follows: its vertices are elements of
$\pi(G)$, and two distinct vertices $r$ and $s$ are adjacent if and only if $rs\in\omega(G)$.
Recall that a subset of vertices of a graph is called a {\it coclique}, if
every two vertices of this subset are non-adjacent. Denote by $t(G)$ the greatest size of a coclique in $GK(G)$.
We refer to a coclique containing $r$ as an $\{r\}$-coclique.
If $r\in\pi(G)$ then $t(r,G)$ is the greatest size of an $\{r\}$-coclique and $\rho(r,G)$ is a set of vertices in some
$\{r\}$-coclique of size $t(r,G)$.

\begin{lemma}\emph{(\cite[Proposition~2]{Vasiliev}, \cite[Theorem~2]{VasGor})}\label{l:Vasiliev theorem}
Let $L$ be a finite nonabelain simple group with $t(L)\geq3$ and $t(2,L)\geq2$, and let $G$
be a finite group isospectral to $L$. Then the following statements hold.
\\$(1)$ There exists a nonabelain simple group $S$ such that
$S\leq\overline{G}=G/K\leq\Aut S$ for the maximal normal soluble subgroup $K$ in $G$.
\\$(2)$ For every coclique $\rho$ of $GK(G)$ containing at least three
elements, at most one prime in $\rho$ divides the product $|K|\cdot|\overline{G}/S|$. In particular,
$t(S)\geq~t(L)-1$.
\\$(3)$ Every prime $r\in\pi(G)$ non-adjacent to $2$ in $GK(G)$
does not divide  $|K|\cdot|\overline{G}/S|$. In particular,
$t(2,S)\geq{t(2,L)}$.
\end{lemma}

The values of $t(S)$ and $t(2,S)$ for all nonabelain simple group $S$ were obtained in \cite{VasVd_criteria,VasVd_cocliques}, in particular, $t(E_7(q))=8$ and
$t(2,E_7(q))\geq3$. Lemma~\ref{l:Vasiliev theorem} shows that the nonabelian composition factor $S$ of a group isospectral to $E_7(q)$ must satisfy $t(S)\geq7$ and $t(2,S)\geq3$. Table~\ref{tab:table1} contains all simple groups $S$ that enjoy such properties. The information in this table is extracted from \cite{VasVd_criteria,VasVd_cocliques}.

\begin{table}%
\centering
{\small
\begin{tabular}{|c|ll|c|@{}c@{}|c|}
\hline
  $S$  & \multicolumn{2}{c|}{Conditions} & $t(2,S)$ & $\rho(2,S)\setminus \{2\}$ & $t(S)$ \\
\hline
 $J_{4}$  & & & 6 & $\{23,29,31,37,43\}$ & 7\\
 $F_{1}$  & none & & $5$ & $\{29,41,59,71\}$ & $11$\\
 $F_{2}$  & & & $3$ & $\{31,47\}$ & $8$\\
\hline
 ${Alt}_n$ & \multicolumn{2}{l|}{$n, n-2$ are prime} & $3$ & $\{n, n-2\}$&\underline{\mathstrut\ \ \ }\\
 $n\ge 47 $& \multicolumn{2}{l|}{$n-1, n-3$ are prime} & $3$ & $\{n-1,n-3\}$ & \\
\hline
 $A^{\tau}_{n-1}(u)$ & \multicolumn{2}{l|}{$2<(u-\tau1)_{\{2\}}=n_{\{2\}}$} & $3$ & $\{r_{n-1}(\tau{u}),r_n(\tau{u})\}$& \\
 $n\ge 13 $& \multicolumn{2}{l|}{$u\equiv0\pmod2$} &$3$ & $\{r_{n-1}(\tau{u}),r_n(\tau{u})\}$&\myraise{$\left[\frac{n+1}{2}\right]$}\\
\hline
 $B_n(u)$, $n\ge 9$  & $u\equiv0\pmod2$  & $n\equiv1\pmod2$   &$3$ & $\{r_n,r_{2n}\}$& $\left[\frac{3n+5}{4}\right]$\\
\hline
 & $u\equiv 5\pmod 8$ & $n\equiv1\pmod2$ & $3$ &$\{r_n,r_{2n-2}\}$& \\
 \cline{2-5}
 \myraise{$D_n(u)$} & $u\equiv0\pmod2$ & $n\equiv0\pmod2$ & $3$ & $\{r_{n-1},r_{2n-2}\}$& $\left[\frac{3n+1}{4}\right]$ \\
 \myraise{$n\ge 9$}& & $n\equiv1\pmod2$ & $3$ &$\{r_{n},r_{2n-2}\}$& \\
\hline
 & $u\equiv 3\pmod 8$ &$n\equiv1\pmod2$& $3$ & $\{r_{2n-2},r_{2n}\}$& \\
 \cline{2-5}
 \myraise{${}^2D_n(u)$}& $u\equiv0\pmod2$ & $n\equiv0\pmod2$ & $4$ & $\{r_{n-1},r_{2n-2}, r_{2n}\}$& $\left[\frac{3n+4}{4}\right]$\\
 \myraise{$n\ge 8$}& & $n\equiv1\pmod2$ & $3$ & $\{r_{2n-2},r_{2n}\}$ & \\
\hline
 &$u\equiv1 \pmod4$&&3&$\{r_{14},r_{18}\}$&\\
 $E_7(u)$&$u\equiv3 \pmod4$&&3&$\{r_{7},r_{9}\}$&8\\
 &$u\equiv0\pmod2$ &&5&$\{r_7,r_9,r_{14},r_{18}\}$&\\
\hline
 $E_8(u)$& none &&5&$\{r_{15},r_{20},r_{24},r_{30}\}$&12\\
\hline
\end{tabular}}
\caption{Simple groups $S$ with $t(S)\geq 7$ and $t(2,S)\geq3$}
\label{tab:table1}
\end{table}

Following \cite{VasVd_cocliques}, by the compact form for the prime graph of a finite simple group $G$ of Lie type over a field of order $q$ and characteristic $p$ we mean a
graph whose vertices are labeled with marks $R_i$ and $p$. The vertex labeled $R_i$ represents a clique of $GK(G)$ such that every vertex in this clique labeled by a prime in $R_i(q)$. An edge
joining $R_i$ and $R_j$ represents the edges of $GK(G)$ that join each vertex in $R_i(q)$ to each vertex in $R_j(q)$. Finally, an edge between $p$ and $R_i$ means that $p$ is adjacent to all primes in~$R_i(q)$. Figure~1 presents the compact form of $GK(E_7(q))$ (see \cite[Figure~4]{VasVd_cocliques}).

\begin{figure}[ht]
\centering
\includegraphics[scale=0.7]{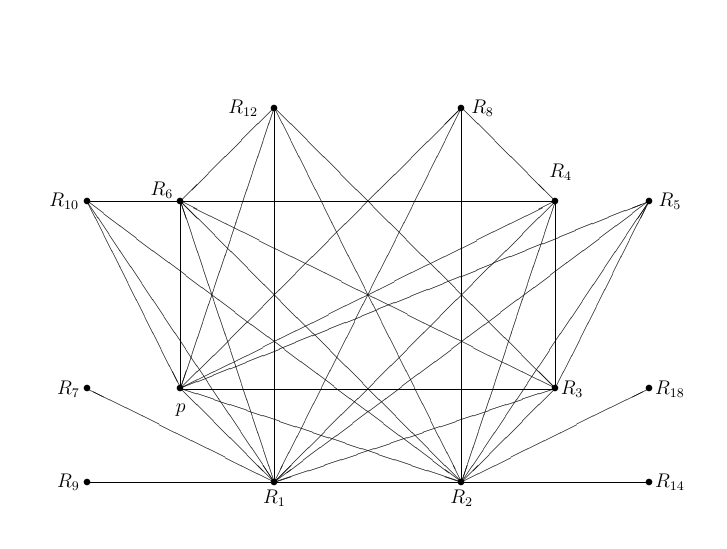}
\caption{Compact form of $GK(E_7(q))$}
\label{fig1}
\end{figure}

\vspace{1\baselineskip}

It is known that the order of any semisimple element of a finite simple group of Lie type divides the order of some maximal torus of this group.
The maximal tori of the universal groups $\overline{E}_7(q)$ were described in \cite{tori}.  Recall that $\overline{E}_7(q)\simeq d.E_7(q)$ is a central extension of the group of order $d=(q-1,2)$ by $E_7(q)$. Table~2 gives a cyclic structure of the maximal tori of~$\overline{E}_7(q)$. In this table, given a nonzero integer $k$, $Z_k$ stands for the cyclic group of order $|k|$, $(Z_k)^m$ means the direct product of $m$ groups isomorphic to $Z_k$, and $\epsilon\in\{+,-\}$.

\begin{table}%
\centering
\renewcommand{\arraystretch}{1.0}
{\small
\begin{tabular}{ | l | l | }
\hline
$(Z_{\epsilon{q}-1})^7$ &  $Z_{\epsilon{q}-1}\times(Z_{\epsilon{q}+1})^2\times Z_{q^4-1}$  \\ \hline
$(Z_{\epsilon{q}-1})^5\times Z_{q^2-1}$ &  $Z_{\epsilon{q}-1}\times Z_{(\epsilon{q}+1)((\epsilon{q})^5-1)}$ \\ \hline
$(Z_{\epsilon{q}-1})^3\times(Z_{q^2-1})^2$ & $Z_{\epsilon{q}-1}\times Z_{q^6-1}$ \\ \hline
$(Z_{\epsilon{q}-1})^4\times Z_{(\epsilon{q})^3-1}$ & $Z_{\epsilon{q}-1}\times Z_{(q^2-1)(q^4+1)}$ \\ \hline
$Z_{\epsilon{q}-1}\times(Z_{q^2-1})^3$ & $Z_{\epsilon{q}-1}\times Z_{(\epsilon{q}-1)(q^2+1)((\epsilon{q})^3+1)}$  \\ \hline
$(Z_{\epsilon{q}-1})^2\times Z_{q^2-1}\times Z_{(\epsilon{q})^3-1}$ & $(Z_{q^2+\epsilon{q}+1})^2\times Z_{(\epsilon{q})^3-1}$ \\ \hline
$(Z_{\epsilon{q}-1})^3\times Z_{q^4-1}$ & $Z_{(\epsilon{q})^3+1}\times Z_{(\epsilon{q})^3-1}\times Z_{\epsilon{q}+1}$ \\ \hline
$Z_{\epsilon{q}-1}\times(Z_{\epsilon{q}+1})^2\times(Z_{q^2-1})^2$ & $Z_{((\epsilon{q})^3-1)(q^4-q^2+1)}$  \\ \hline
$Z_{\epsilon{q}-1}\times Z_{q^2-1}\times Z_{(\epsilon{q}+1)((\epsilon{q})^3-1)}$ & $Z_{(\epsilon{q}-1)(q^6+\epsilon{q}^3+1)}$ \\ \hline
$Z_{\epsilon{q}-1}\times (Z_{(\epsilon{q})^3-1})^2$ & $Z_{q^2-\epsilon{q}+1}\times Z_{(\epsilon{q}-1)(q^4+q^2+1)}$ \\  \hline
$Z_{\epsilon{q}-1}\times Z_{q^2-1}\times Z_{q^4-1}$ & $Z_{(\epsilon{q})^3-1}\times Z_{q^4-1}$  \\ \hline
$(Z_{\epsilon{q}-1})^2\times Z_{(\epsilon{q})^5-1}$ & $Z_{((\epsilon{q})^5-1)(q^2+\epsilon{q}+1)}$ \\ \hline
$Z_{\epsilon{q}-1}\times Z_{q^2-1}\times Z_{(\epsilon{q}-1)((\epsilon{q})^3+1)}$ & $Z_{(\epsilon{q}-1)(q^2+1)}\times Z_{q^2-1}\times Z_{q^2+1}$  \\ \hline
$Z_{\epsilon{q}-1}\times(Z_{(\epsilon{q}-1)(q^2+1)})^2$ &  $Z_{(\epsilon{q})^7-1}$ \\ \hline
$Z_{(\epsilon{q})^3-1}\times Z_{(\epsilon{q}+1)((\epsilon{q})^3-1)}$ & $Z_{q^4+1}\times Z_{(\epsilon{q}-1)(q^2+1)}$ \\ \hline
\end{tabular}}
\caption{Maximal tori of $(2,q-1).E_7(q)$}
\label{tab:table2}
\end{table}

\begin{lemma}\label{2coclique} Let $S$ be a finite simple group of Lie type over a field of order $u$ from Table~\emph{\ref{tab:table1}}.
Suppose that $\{2,r_{i_1}(u),r_{i_2}(u),...,r_{i_m}(u) \}$ is a $\{2\}$-coclique in $GK(S)$ of the greatest size
and define $M=\{i_1,i_2,...,i_m\}$. If $r_j(u)$ is non-adjacent to $2$ in $GK(S)$, then $j\in{M}$. Moreover, $k_j(u)$ is the maximal w.r.t. divisibility $R_j(u)$-number in $\omega(S)$.
\end{lemma}
\begin{prf}
If $u$ is even, then the first assertion follows from \cite[Proposition 6.4]{VasVd_criteria}. In the odd case, it holds by \cite[Proposition 6.7]{VasVd_criteria}. The cyclic structure of tori in groups of Lie type from Table~1 is known (see \cite{ButGr} for classical groups and \cite{tori} for exceptional groups). This gives sufficient information to verify the last statement of the lemma.
\end{prf}

\begin{lemma}\label{divisible}
Let $G$ be a finite group isospectral to $E_7(q)$, $q$ odd, and $S\leq\overline{G}=G/K\leq\Aut S$, where $K$
is the soluble radical of $G$ and $S$ is a simple group of Lie type over a field of order $u$ and characteristic~$v$. Choose $\varepsilon\in\{+,-\}$ such that $q\equiv{-\varepsilon}1\pmod4$.
Then there exist integers $i_1$, $i_2$, $m_1(S),m_2(S)\in\mu(S)$ such that the following chains of divisibilities hold
$$(q^7-\varepsilon1)/2~\vdots~{m_1(S)}~\vdots~{k_{i_1}(u)}~\vdots~{k_7(\varepsilon{q})};$$
$$(q-\varepsilon1)\cdot(q^6+\varepsilon{q^3}+1)/2~\vdots~{m_2(S)}~\vdots~{k_{i_2}(u)}~\vdots~{k_9(\varepsilon{q})}.$$

\end{lemma}
\begin{prf} It follows from Lemma~\ref{l:Vasiliev theorem} that $S$ is among the groups from Table~1.
Since $q\equiv{-\varepsilon}1\pmod4$, a set $\{2,r_7(\varepsilon{q}),r_9(\varepsilon{q})\}$ is a coclique in $GK(L)$ (see Table~\ref{tab:table1} and take into account that $k_{2n}(q)=k_n(-q)$ for odd $n$). Therefore, Lemma~\ref{l:Vasiliev theorem} yields that $r_i(\varepsilon{q})\in\pi(S)$ and $(r_i(\varepsilon{q}),|\overline{G}/S||K|)=1$ for $i=7,9$. Note that there exists a coclique of size 8 in $GK(L)$ which contains $r_7(\varepsilon{q})$ and $r_9(\varepsilon{q})$, hence
$t(r_i(\varepsilon{q}),S)\geq7$ for $i=7,9$ due to Lemma~\ref{l:Vasiliev theorem}. It follows from \cite[Tables~4 and~5]{VasVd_criteria} that $t(v,S)\leq5$.
Therefore $v\not\in R_7(\varepsilon{q})\cup R_9(\varepsilon{q})$. Thus there exist indices $i_1$, $i_2$ such that $r_7(\varepsilon{q})\in{R}_{i_1}(u)$,
$r_9(\varepsilon{q})\in{R}_{i_2}(u)$. Moreover, $R_7(\varepsilon{q})\subseteq{R_{i_1}(u)}$ and $R_9(\varepsilon{q})\subseteq{R_{i_2}}(u)$. Indeed, if $r$ and $s$ are two distinct primes in $R_7(\varepsilon{q})$ (or  $R_9(\varepsilon{q})$), then $r$ and $s$ are adjacent in $GK(L)$, so they are adjacent in $GK(S)$ by preceding arguments. On the other hand, $r$ and $s$ are non-adjacent to $2$ in $GK(S)$, so Lemma~\ref{2coclique} implies that $e(r,u)=e(s,u)$.

It follows from Lemma~\ref{2coclique} that $k_7(\varepsilon{q})\in\omega(L)$ and $k_{i_1}(u)\in\omega(S)$. Moreover, $(k_7(\varepsilon{q}),|\overline{G}/S||K|)=1$ and $k_{i_1}(u)$ is the maximal w.r.t. divisibility $R_{i_1}(u)$-number in $\omega(S)$, so $k_7(\varepsilon{q})~|~k_{i_1}(u)$. Similarly, $k_9(\varepsilon{q})~|~k_{i_2}(u)$. Obviously, there exist $m_1(S)$ and $m_2(S)$ in $\mu(S)$ such that $k_{i_1}(u)~|~m_1(S)$ and $k_{i_2}(u)~|~m_2(S)$. Since $\omega(S)\subseteq\omega(G)$, the numbers $m_1(S)$, $m_2(S)$ lie in $\omega(G)$. Primes $r_7(\varepsilon{q})$ and $r_9(\varepsilon{q})$ are non-adjacent to~$p$ in $GK(G)$ (see Figure \ref{fig1}), so elements of orders $m_1(S)$ and $m_2(S)$ are semisimple in~$L$. Hence $m_1(S)$ and $m_2(S)$ divide orders of some maximal tori of $L$. By Table \ref{tab:table2}, $k_7(\varepsilon{q})$ divides only the integer $(q^7-\varepsilon{1})/2$ among the orders of maximal tori, so $m_1(S)~|~(q^7-\varepsilon{1})/2$. By the same reason, $m_2(S)~|~(q-\varepsilon1)(q^6+\varepsilon{q}^3+1)/2$.
\end{prf}

\begin{lemma}\label{unipotent}
Suppose that $L$ is a finite simple group of Lie type over a field of characteristic $p$ and $\exp_p(L)$ is the exponent of a Sylow $p$-subgroup of $L$. Then $\exp_p(L)=\min\{p^\alpha~|~p^\alpha>\operatorname{ht}(L)\}$, where $\operatorname{ht}(L)$ is the height of the highest root in the root system of $L$. In particular, if $L$ is of type $E_7$ then $\exp_p(L)=\min\{p^\alpha~|~p^\alpha>17\}$.
\end{lemma}

\begin{prf} This follows from \cite[Corollary 0.5]{Test}.
\end{prf}

\begin{lemma}\label{l:reduction}\emph{\cite[Lemma~2.3]{cover}}
Let $A$ and $B$ be finite groups. The following are equivalent.
\newline$(1)$ $\omega(H)\not\subseteq\omega(B)$ for any proper cover $H$ of $A$.
\newline$(2)$ $\omega(H)\not\subseteq\omega(B)$ for any split extension $H=K:A$, where $K$ is a nontrivial elementary abelian group.
\end{lemma}

\begin{lemma}\emph{{\cite[Lemma~1.5]{VasGrSt}}}\label{l:action2} Let $G$ be a finite group, $K$ be a normal subgroup of $G$, and
$r\in\pi(K)$. Suppose that the factor group $G/K$ has a section isomorphic to a non-cyclic abelian $p$-group for some odd prime $p$ distinct from~$r$. Then
$rp\in\omega(G)$.
\end{lemma}

If $G$ is a group, $g\in G$, and $V$ is a finite-dimensional $G$-module, then $\deg_V(g)$ stands for the degree of the minimal polynomial of $g$ on $V$. The next assertion is well-known.

\begin{lemma}\label{MinimalPolynomial}
Suppose that $G$ is a finite group, $V$ is a finite-dimensional $G$-module over a field of positive characteristic $r$ and $H=V\leftthreetimes{G}$ is a natural semidirect product. The orders of elements of a coset $Vg$ in $H$ coincide with the order of $g$ in $G$ if and only if the minimal polynomial of $g$ on $V$ divides $(x^{|g|}-1)/(x-1)$. In particular, if $\deg_V(g)=|g|$, then $Vg$ contains an element of order $r|g|$.
\end{lemma}

\begin{lemma}\label{Structure2Dn}
Let $S={}^2D_n(u)$, $n\geq8$, and $u=v^m$ for a prime $v$. Then the following statements hold.
\newline$(1)$ If $r\in\pi(S)$ and $r$ does not divide the order of any proper parabolic subgroup of $S$, then $e(r,u)=2n$.
\newline$(2)$ If $V$ is a finite-dimensional $S$-module over a field of characteristic $r\neq v$ and $g$ is an element of prime order $s\neq v$ which lies in some proper parabolic subgroup of $S$, then $\deg_V(g)=s$. In particular, the natural semidirect product $V\leftthreetimes{S}$ contains an element of order $rs$.
\end{lemma}

\begin{prf} (1) The orders and structure of parabolic subgroups of finite classical groups are well-known (see, e.g., \cite[Proposition~4.20]{KL}).

(2) We may assume that $V$ is absolutely irreducible. Then the first statement follows from the main theorem of \cite{DMZal}. The second assertion follows from the first one and Lemma~\ref{MinimalPolynomial}.

 \end{prf}

\sect{Proof of the theorem: a nonabelain composition factor}

Let $L=E_7(q)$, $q=p^m$, $p$ a prime, and $G$ be a finite group with $\omega(G)=\omega(L)$. Since
the quasirecognizability of $E_7(2)$ and $E_7(3)$ was proved in \cite{AlKon03}, further we assume that $q\geq4$.
Note that $GK(G)=GK(L)$, so $t(2,G)\geq3$ and $t(G)=8$. By Lemma \ref{l:Vasiliev theorem},
there is a nonabelain simple group $S$ such that
$S\leq{G/K}\leq\Aut S$ for the maximal normal soluble subgroup $K$ of $G$, 
$t(S)\geq7$, and $t(2,S)\geq3$. Thus $S$ is one of the groups given in Table~\ref{tab:table1}.
We consider every case separately and show that $S\simeq{L}$.

\begin{lemma}\label{l:sporadic}$S$ is not isomorphic to a sporadic group or the Tits group.\end{lemma}
\begin{prf}
It is a direct consequence of \cite[Lemma~7]{sporadic}.
\end{prf}

\begin{lemma}\label{even case} If $p=2$, then $S\simeq{L}$.
\end{lemma}
\begin{prf}
In this case $t(2,L)=5$, so it follows from Lemma~\ref{l:Vasiliev theorem} that
$t(2,S)\geq5$. Using Table \ref{tab:table1} we find that $S$ can be isomorphic to either
$E_7(2^{k})$, or $E_8(u)$, or sporadic groups $J_4$ and $F_1$. By Lemma \ref{l:sporadic}, only
the cases $E_7(2^{k})$ and $E_8(u)$ are possible.

Let $S\simeq E_8(u)$ and $u=v^{k}$. Applying Lemma~\ref{unipotent} for $L$ and its subsystem subgroups, we derive $32\in\mu(L)$.
On the other hand, in the group $E_8(u)$ there is an element of order $32s$, where $s$ is an odd prime. Indeed, if $u$ is even then $S$ contains elements of order $32(u\pm1)$ \cite{Deriziotis}, and if $v$ is odd then $u^8-1\in\omega(S)$ \cite{tori}; a contradiction.

Let $S\simeq E_7(u)$ and $u=2^{k}$. Note that if $i=7,9,14,18$ then $r_i(q)$ is non-adjacent to 2 in $GK(L)$ (see for example Figure \ref{fig1}),
so Lemma \ref{l:Vasiliev theorem} implies that $r_i(q)\in\pi(S)$ and $(r_i(q),|K||\overline{G}/S|)=1$.
For $i=7,9,14,18$ choose a primitive prime divisor $r_i\in{R}_i(q)$ such that
$e(r_i,2)=im$. Put $e_i=e(r_i,u)$. Then $r_i$ divides $u^{e_i}-1=2^{e_ik}-1$. Therefore $im$ divides
$e_ik$. Suppose $e_{18}k>18m$. Since $k_{e_{18}}(u)$ divides $|S|$, a prime
$r$ with $e(r,2)=e_{18}k$ lies in $\omega(S)$. However, $e(r,q)>18$, so $r\notin\omega(L)$; a contradiction. Thus
$e_{18}k=18m$. If $e_{14}k>14m$ then
$e_{14}k\geq2\cdot14m>18m$ which is impossible by the
same reason. Similarly, $e_9k$ is equal to $9m$ or $18m$.
However, $e_9k\neq e_{18}k=18m$, so we have
$e_9k=9m$. Finally, we derive that $e_7k=7m$. In
particular, $e_{18}>e_{14}>e_9>e_7$. On the other hand,
$e_{18},e_{14},e_9,e_7\in\{18,14,9,7\}$. Thus, $e_{18}=18$ and
$k=m$. The lemma is proved.
\end{prf}

From this moment, we may assume that $q$ is odd. Fix $\varepsilon\in\{+,-\}$ such that $q\equiv{-\varepsilon}1\pmod4$. Then $\{2,r_7(\varepsilon{q}),r_9(\varepsilon{q})\}$ is a coclique in $GK(G)$ (see Lemma~\ref{divisible}).

\begin{lemma} $S\not\simeq{Alt}_n$. \end{lemma}
\begin{prf}
Assume that $S\simeq{Alt}_{n}$. Table \ref{tab:table1} implies $t(2,S)=3$. Let $r$ and $r+2$ be primes in
$\{n-3,n-2,n-1,n\}$. Using Lemma \ref{l:Vasiliev theorem} we see that $k_{9}(\varepsilon{q})$ divides $r$ or $r+2$. Since $q>3$,
the inequality $q^6\geq16q^4$ holds. It follows that $k_{9}(\varepsilon{q})=\frac{q^6+\varepsilon{q}^3+1}{(q-\varepsilon1,3)}>q^4+2$. Hence $q^4\in\omega(Alt_n)$.
On the other hand, $q^3>17$ and, by Lemma~\ref{unipotent}, $q^4\notin\omega(L)$; a contradiction.
\end{prf}

Thus, we may assume that $S$ is a group of Lie type. Suppose further that $S$ is defined over a field of order $u$, where $u=v^k$ for a prime $v$ and a positive integer~$k$.

\begin{lemma} $S\not\simeq A^{\tau}_{n-1}(u)$.
\end{lemma}

\begin{prf} Assume the opposite. It follows from \cite[Table 3]{VasVd_cocliques} that
$t(S)=[\frac{n+1}{2}]$. Then $t(S)\geq7$ provides $n\geq13$. Moreover, Lemma \ref{l:Vasiliev theorem} and \cite[Tables~4 and~6]{VasVd_criteria}
imply that $t(2,S)=3$ and $\rho(2,S)=\{2,r_{n-1}(\tau{u}),r_{n}(\tau{u})\}$.
One of the numbers $n-1$ or $n$ must be even. Let $n-1$ be even. By Lemma \ref{divisible},
there exists $i\in\{7,9\}$ such that $r_{i}(\varepsilon{q})\in R_{n-1}(\tau{u})$.
Let $m_7(L)=(q^7-\varepsilon1)/2$,
$m_9(L)=(q^6+\varepsilon{q}^3+1)(q-\varepsilon1)/2$.
By Lemma \ref{divisible}, $k_{i}(\varepsilon{q})$ divides $k_{n-1}(\tau{u})$, $k_{n-1}(\tau{u})$
divides $m_1(S)$, and $m_1(S)$ divides $m_i(L)$. Since $n-1$ is even, the equality $(\tau{u})^{n-1}-1=((\tau{u})^{(n-1)/2}-1)((\tau{u})^{(n-1)/2}+1)$ holds,
where $(n-1)/2$ is an integer. Now the definition of $k_{n-1}(\tau{u})$ implies that $k_{n-1}(\tau{u})\leq|(\tau{u})^{(n-1)/2}+1|$.
On the other hand, $m_1(S)=\frac{u^{n-1}-1}{(n,\tau{u}-1)}$ \cite[Corollary 3]{LinSpec}. Furthermore, $\frac{u^{n-1}-1}{(n,\tau{u}-1)}\geq\left|\frac{u^{n-1}-1}{\tau{u}-1}\right|=
\left|\frac{(\tau{u})^{(n-1)/2}-1}{\tau{u}-1}((\tau{u})^{(n-1)/2}+1)\right|>|u^{(n-5)/2}(\tau{u}+1)((\tau{u})^{(n-1)/2}+1)|$. However,
$|u^{(n-5)/2}(\tau{u}+1)|>|(\tau{u})^{(n-1)/2}+1|^{1/2}$ due to $n\geq13$. Therefore,
$m_1(S)>|(\tau{u})^{(n-1)/2}+1|^{3/2}>k_{n-1}(\tau{u})^{3/2}$. It follows that
$k_{i}(\varepsilon{q})^{3/2}\leq k_{n-1}^{3/2}(\tau{u})\leq m_1(S)\leq m_i(L)$.
One may easily verify that it is impossible for $i=7,9$.
If $n$ is even, then $k_{n}(\tau{u})\leq |(\tau{u})^{n/2}+1|$ and the same argument gives us a contradiction.

\end{prf}

\begin{lemma}\label{BnCn} $S\not\simeq B_n(u)$ and $S\not\simeq C_n(u)$. \end{lemma}
\begin{prf}
Let $S\simeq B_{n}(u)$ or $S\simeq C_{n}(u)$. Then $u$ is even, $B_n(u)\simeq C_n(u)$, and $n\geq9$. Note that $t(r,L)\geq3$ for every $r\in\pi(L)$ (see Figure~1), in particular, $t(3,L)\geq3$. In fact, if $3$ divides
$q+\epsilon1$ then $3$ is non-adjacent to $r_7(\epsilon{q})$ and $r_9(\epsilon{q})$; while if $p=3$ then it is non-adjacent to every $r_{i}(\pm{q})$, where $i=7,9$.
An adjacency criterion for the prime graph of groups $B_{n}(u)$
implies that $t(3,S)=2$ (see \cite[Proposition~3.1]{VasVd_criteria} and \cite[Proposition~2.4]{VasVd_cocliques}). It follows that one of the primes from $\rho(3,L)\setminus\{3\}$, say $r$, should be coprime to $|S|$. Observe that a Sylow $3$-subgroup of $S$ is non-cyclic due to $n\geq9$. So if $r\in\pi(K)$, we derive a contradiction by Lemma~\ref{l:action2}. Therefore, one of the numbers from $\{k_7(\pm{q}), k_9(\pm{q})\}$ divides $|\Out{S}|$. Since $u$ is even, we have $|\Out{S}|=k$. Therefore, $k\geq\min\{k_7(\pm{q}),k_9(\pm{q})\}$. However, $k_7(\pm{q})\geq(q^6-q^5+q^4-q^3+q^2-q+1)/7\geq(5q^5-q^5+q^4-q^3+q^2-q+1)/7>q^5/2$
and $k_9(\pm{q})\geq(q^6-q^3+1)/3\geq(5q^5-q^3+1)/3>q^5$. So $2k>q^5$. The inequality $n\geq9$ yields that $u^4-1,u^4+1\in\omega(S)\subseteq\omega(L)$.
At least one of these numbers is not divisible by $p$ and so it is the order of a semisimple element of $L$.
On the other hand, $u^4-1=2^{4k}-1=(2^{2k})^2-1\geq(2k)^2-1>q^{10}-1$. However, $q^{10}-1$ is greater than every number in Table 2; a contradiction.
\end{prf}

\begin{lemma}\label{Dn} $S\not\simeq D_n^{\tau}(u)$, where $\tau\in\{+,-\}$.
\end{lemma}
\begin{prf}

(1) Assume the opposite. By \cite[Proposition~3.1]{VasVd_criteria} and \cite[Proposition~2.5]{VasVd_cocliques}, we have $t(3,S)<3$ whenever $v\neq3$. Therefore, if $v\neq3$, then one of the numbers in $\{k_7(\pm{q}), k_9(\pm{q})\}$ divides $|\Out{S}|$, and we derive a contradiction as we did it in Lemma~\ref{BnCn}. It follows that $v=3$.
If $S\simeq D_{n}(u)$, then $t(2,S)=3$ implies $u\equiv 5\pmod8$, which is obviously impossible for $v=3$. Suppose that $S\simeq{}^2D_{n}(u)$.  Then $t(2,S)=3$
if and only if $n$ is odd and $u\equiv3\pmod8$. It follows that $S\simeq {}^2D_{n}(u)$, $n$ is odd, and $u=3^{k}$,
$k$ is odd. Note that in this case $\Out{S}$ is a group of order $8k$ and $\rho(2,S)=\{2,r_{2n-2}(u),r_{2n}(u)\}$.

(2) Recall that we fix $\varepsilon\in\{+,-\}$ such that $\rho(2,L)=\{2,r_7(\varepsilon{q}),r_9(\varepsilon{q})\}$. Put $\sigma=\{r_7(q),r_9(q),r_7(-q),r_9(-q)\}$. Then $\sigma$ is a coclique in $GK(L)$, so at most one prime from $\sigma$ can divide the product $|K||\overline{G}/S|$, furthermore, such a prime is adjacent to~$2$. Among the remaining three numbers only one, say $t$, belongs to $R_{2n}(u)$, moreover, $t$ is non-adjacent to~$2$. Since every $r\in\pi(L)$ is non-adjacent to either $r_7(\varepsilon{q})$ and $r_9(\varepsilon{q})$ or $r_7(-\varepsilon{q})$ and $r_9(-\varepsilon{q})$ (see Figure~1), it follows that for every $r\in\pi(L)$ there exists a prime $s$ from $\sigma$ such that $s$ and $r$ are non-adjacent in $GK(L)$, $s$ is coprime to $|K||\overline{G}/S|$, and $s\not\in R_{2n}(u)$, in particular, by Lemma~\ref{Structure2Dn}(1), $s$ divides the order of some proper parabolic subgroup of~$S$.

(3) Assume that $u=3^k>3$. Since $k$ is odd, there exists an odd prime $r$ lying in ${R}_1(u)$. It follows from \cite[Proposition~2.5]{VasVd_cocliques} that $t(r,S)=2$ and $r$ is non-adjacent to $t$ in $GK(S)$ if and only if $t\in R_{2n}(u)$. This contradicts~(2).

(4) Thus, $u=3$. Suppose that the soluble radical $K$ of $G$ is non-trivial. We claim that $\omega(G)\not\subseteq\omega(L)$ in this case. If $H$ is the preimage of $S$ in $G$, then $H$ is a proper cover of $S$. By Lemma~\ref{l:reduction}, in order to prove that $\omega(H)\not\subseteq\omega(L)$ it is sufficient to show that $\omega(V:S)\not\subseteq\omega(L)$, where $V$ is an elementary abelian $r$-subgroup for some prime~$r$. So we may assume that $K=V$ and $G=V:S$. If $C_G(K)\nleq K$, then $G=C_G(K)$ due to simplicity of $S$, and $r$ is adjacent to every prime in $GK(G)$, which contradicts~(2). So $S$ acts faithfully on $K$. Choose for the prime $r$ a prime $s$ as in~(2). If $r\neq 3$, then $rs\in\omega(G)$ due to Lemma~\ref{Structure2Dn}(2). By \cite{GT}, the group $S={}^2D_{n}(3)$ for odd $n$ is
unisingular, that is every its semisimple element has a non-trivial fixed point on every abelian $3$-subgroup $K$ with $S$-action, so if $r=3$ then $rs\in\omega(G)$ as well. On the other hand, $rs\not\in\omega(L)$ by the choice of~$s$. Thus,  $K$ must be trivial. Since $u=3$, the order of $\Out S$ is equal to~$8$. So the inequality $t(L)=8>t(2,L)=t(2,S)$ implies that $t(S)=t(L)=8$. However, $n$ is odd and $n\geq9$, hence $t(S)=\left[\frac{3n+4}{4}\right]\geq9$ for $n\geq11$ and $t(S)=7$ for $n=9$; a contradiction. This completes the proof.
\end{prf}

Thus, $S$ should be a finite simple exceptional group of Lie type.

\sect{Completion of the proof}

By preceding arguments, we have that $S\leq\overline{G}=G/K\leq\Aut S$, where $K$~is the soluble radical of~$G$ and either $S\simeq{E_7(u)}$ or $S\simeq{E_8(u)}$. It appears that the case of $S\simeq{E_8(u)}$ requires a careful study, so it is convenient to have a structure of $GK(E_8(u))$. Figure~2 presented below is taken from \cite[Figure~5]{VasVd_cocliques} and gives a compact form of the prime graph of $E_8(u)$ (the definition of the compact form is formulated before Figure~\ref{fig1}). Observe that the vector from 5 to $R_4$ and the dotted edge $(5,R_{20})$ mean that $R_4$ and $R_{20}$ are not connected, but if $5\in{R_4}$
(i.e., $q^2\equiv-1\pmod5$), then there exists an edge between 5 and~$R_{20}$.

\begin{figure}[ht]
\centering
\includegraphics[scale=0.7]{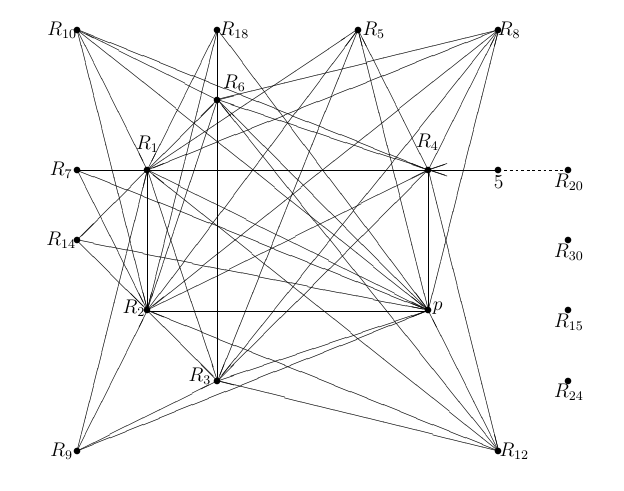}
\caption{Compact form for $GK(E_8(u))$}
\label{fig2}
\end{figure}

The maximal tori of the group $E_8(q)$ were described in \cite{tori}. Table~3 gives a cyclic structure for these
tori. Here we use the same notation as in Table~\ref{tab:table2}.

\begin{table}%
\centering
\renewcommand{\arraystretch}{0.9}
{\small
\begin{tabular}{ | l | l | }
\hline
$(Z_{\epsilon{q}-1})^8$ &  $Z_{(\epsilon{q}+1)((\epsilon{q})^3-1)}\times Z_{q^4-1}$ \\  \hline
$(Z_{\epsilon{q}-1})^6\times Z_{q^2-1}$ & $Z_{q^4-1}\times Z_{q^4-1}$ \\  \hline
$(Z_{\epsilon{q}-1})^4\times(Z_{q^2-1})^2$ & $(Z_{q^2-1})^2\times(Z_{q^2+1})^2$ \\  \hline
$(Z_{\epsilon{q}-1})^5\times Z_{(\epsilon{q})^3-1}$ &  $Z_{q^2-1}\times Z_{(\epsilon{q}+1)((\epsilon{q})^5-1)}$ \\  \hline
$(Z_{\epsilon{q}-1})^2\times(Z_{q^2-1})^3$ &  $Z_{q^2-1}\times Z_{q^6-1}$ \\  \hline
$(Z_{\epsilon{q}-1})^3\times Z_{q^2-1}\times Z_{(\epsilon{q})^3-1}$ &  $Z_{q^2-1}\times Z_{q^6-1}$ \\  \hline
$(Z_{\epsilon{q}-1})^4\times Z_{q^4-1}$ & $Z_{(\epsilon{q}-1)(q^2+1)}\times Z_{(q^2+1)((\epsilon{q})^3-1)}$ \\ \hline
$(Z_{q-1})^2\times(Z_{q+1})^2\times(Z_{q^2-1})^2$ &  $Z_{q^2-1}\times Z_{(q^2-1)(q^4+1)}$ \\ \hline
$(Z_{\epsilon{q}-1})^2\times Z_{q^2-1}\times Z_{(\epsilon{q}+1)((\epsilon{q})^3-1)}$ & $(Z_{q^2+\epsilon{q}+1})^2\times Z_{(\epsilon{q}+1)((\epsilon{q})^3-1)}$ \\ \hline 
$(Z_{\epsilon{q}-1})^2\times(Z_{(\epsilon{q})^3-1})^2$ & $Z_{(\epsilon{q}+1)(q^2+\epsilon{q}+1)((\epsilon{q})^5-1)}$ \\ \hline
$(Z_{\epsilon{q}-1})^2\times Z_{q^2-1}\times Z_{q^4-1}$ &  $Z_{(\epsilon{q}+1)(q^2+1)((\epsilon{q})^5-1)}$ \\ \hline
$(Z_{\epsilon{q}-1})^3\times Z_{(\epsilon{q})^5-1}$ &  $Z_{(\epsilon{q}+1)((\epsilon{q})^7-1)}$ \\ \hline
$(Z_{\epsilon{q}-1})^2\times Z_{q^2-1}\times Z_{(\epsilon{q}-1)((\epsilon{q})^3+1)}$ & $Z_{q^8-1}$ \\ \hline 
$(Z_{\epsilon{q}-1})^2\times (Z_{(\epsilon{q}-1)(q^2+1)})^2$ & $Z_{q^2-1}\times Z_{q^2+1}\times Z_{q^4+1}$ \\ \hline
$Z_{\epsilon{q}-1}\times Z_{(\epsilon{q})^3-1}\times Z_{(\epsilon{q}+1)((\epsilon{q})^3-1)}$ & $(Z_{q^2+1})^2\times Z_{q^4-1}$  \\ \hline
$(Z_{q-1})^2\times(Z_{q+1})^2\times Z_{q^4-1}$ & $Z_{(\epsilon{q}+1)((\epsilon{q})^3-1)(q^4+1)}$ \\ \hline
$(Z_{\epsilon{q}-1})^2\times Z_{(\epsilon{q}+1)((\epsilon{q})^5-1)}$ & $Z_{(q^2+1)(q^6-1)}$   \\ \hline
$(Z_{\epsilon{q}-1})^2\times Z_{q^6-1}$ & $Z_{(q^2-1)(q^2+\epsilon{q}+1)(q^4-q^2+1)}$  \\ \hline
$(Z_{\epsilon{q}-1})^2\times Z_{(q^2-1)(q^4+1)}$ & $Z_{(q^2-1)(q^6+\epsilon{q}^3+1)}$ \\ \hline
$(Z_{\epsilon{q}-1})^2\times Z_{(\epsilon{q}-1)(q^2+1)((\epsilon{q})^3+1)}$ & $(Z_{q^2-\epsilon{q}+1})^2\times Z_{(\epsilon{q}+1)((\epsilon{q})^3-1)}$ \\ \hline
$Z_{\epsilon{q}-1}\times(Z_{q^2+\epsilon{q}+1})^2\times Z_{(\epsilon{q})^3-1}$ & $Z_{(q^2-1)(q^6+1)}$ \\ \hline
$Z_{q-1}\times Z_{q^3+1}\times Z_{q^3-1}\times Z_{q+1}$ & $(Z_{q^2+\epsilon{q}+1})^4$  \\ \hline
$Z_{\epsilon{q}-1}\times Z_{((\epsilon{q})^3-1)(q^4-q^2+1)}$ & $(Z_{q^4+\epsilon{q}^3+q^2+\epsilon{q}+1})^2$ \\ \hline
$Z_{\epsilon{q}-1}\times Z_{(\epsilon{q}-1)(q^6+\epsilon{q}^3+1)}$ & $Z_{q^2+\epsilon{q}+1}\times Z_{q^6+\epsilon{q}^3+1}$ \\ \hline
$Z_{\epsilon{q}-1}\times Z_{q^2-\epsilon{q}+1}\times Z_{(\epsilon{q}-1)(q^4+q^2+1)}$ & $(Z_{q^2+1})^4$ \\ \hline
$Z_{\epsilon{q}-1}\times Z_{(\epsilon{q})^3-1}\times Z_{q^4-1}$ & $Z_{q^2+1}\times Z_{q^6+1}$ \\ \hline
$Z_{\epsilon{q}-1}\times Z_{((\epsilon{q})^5-1)(q^2+\epsilon{q}+1)}$ & $(Z_{q^4+1})^2$  \\ \hline
$Z_{q^2-1}\times(Z_{(q^2+1)(\epsilon{q}-1)})^2$ & $Z_{(q^4-q^2+1)(q^2+\epsilon{q}+1)}\times Z_{q^2+\epsilon{q}+1}$ \\ \hline
$Z_{\epsilon{q}-1}\times Z_{(\epsilon{q})^7-1}$ &  $Z_{q^4+q^2+1}\times Z_{q^2+q+1}\times Z_{q^2-q+1}$  \\ \hline
$Z_{(\epsilon{q}-1)(q^4+1)}\times Z_{(\epsilon{q}-1)(q^2+1)}$ & $Z_{q^8+\epsilon{q}^7-\epsilon{q}^5-q^4-\epsilon{q}^3+\epsilon{q}+1}$ \\ \hline
$(Z_{q^2-1})^4$ & $Z_{q^8-q^4+1}$ \\ \hline
$(Z_{q^2-1})^2\times Z_{(\epsilon{q}+1)((\epsilon{q})^3-1)}$ & $Z_{q^8-q^6+q^4-q^2+1}$ \\ \hline
$(Z_{q^2-1})^2\times Z_{q^4-1}$ & $(Z_{q^4-q^2+1})^2$ \\ \hline
$Z_{(\epsilon{q}+1)((\epsilon{q})^3-1)}\times Z_{(\epsilon{q}+1)((\epsilon{q})^3-1)}$ &  \\ \hline
\end{tabular}}
\caption{Maximal tori of $E_8(q)$}
\label{tab:table3}
\end{table}

\begin{lemma}\label{R-inclusions}
Let $S\simeq E_7(u)$ or $S\simeq E_8(u)$. Suppose that $I$ and $J$ are subsets of positive integers such that $\bigcup\limits_{i\in I}R_{i}(u)\subseteq\pi(S)$, $\bigcup\limits_{j\in J}R_{j}(q)\subseteq\pi(L)$, and $\bigcup\limits_{i\in I}R_{i}(u)\subseteq\bigcup\limits_{j\in J}R_{j}(q)$. Then $\prod\limits_{i\in I}k_{i}(u)$ divides $d\cdot\prod\limits_{j\in J}k_{j}(q)$, where $d=35$ provided $I\cap\{3,4,6\}\neq\varnothing$,
and $d=1$ otherwise.
\end{lemma}
\begin{prf}
Let $r\in\bigcup\limits_{i\in I}R_{i}(u)$. Then $r\in\bigcup\limits_{j\in J}R_{j}(q)$. So there exist integers $a\in I$ and $b\in J$ such that $r\in{R}_{a}(u)$ and $r\in{R}_{b}(q)$.  Set $r^\alpha=|k_a(u)|_{\{r\}}$ and $r^\beta=|k_b(q)|_{\{r\}}$. In order to prove the lemma it is sufficient to show that $\alpha\leq\beta+1$ if $a\in\{3,4,6\}$ and $r\in\{5,7\}$, and $\alpha\leq\beta$ in all other cases. Assume that the lemma is wrong. Then $\alpha>\beta$, in particular, $\alpha\geq2$. The cyclic structure of maximal tori in simple groups of types $E_7$ and $E_8$ (see Tables~2 and~3) implies that $k_a(u)$ and, consequently, $r^\alpha$ lie in $\omega(S)$, so $r^\alpha\in\omega(L)$. Let $c$ be a least positive integer such that $r^\alpha$ divides $q^c-1$. Then $c\leq18$ (see Table~2). On the other hand, since $e(r,q)=b$, it follows that $c=bf$, where $f$ is a positive integer, and $f$ is greater than~$1$ due to $\alpha>\beta$.

Suppose firstly that $r$ is odd. Observe that $a$ and $b$ divide $r-1$ by Fermat's little theorem. Lemma \ref{r-part} yields that $(q^c-1)_{\{r\}}=(q^b-1)_{\{r\}}\cdot{f}_{\{r\}}$, so $r$ divides~$f$. It is easy to verify using Table~2, that any prime divisor of $c$ does not exceed~$7$. Therefore, $r\leq7$. Suppose that either $r=5$ or $r=7$. By Fermat's little theorem, $a\in\{1,2,3,4,6\}$. If $a=1$ then $r^{\alpha+1}~|~u^{r}-1$, while $a=2$ implies $r^{\alpha+1}~|~u^{r}+1$. Since $(u^{r}-1)/(u-1,2)\in\omega(S)$ and $(u^{r}+1)/(u-1,2)\in\omega(S)$ (see Table~2), it follows that $r^{\alpha+1}\in\omega(L)$. However, if $g$ is the least positive integer such that $r^{\alpha+1}$ divides $q^g-1$, then $g\geq\beta{r}^2$ by Lemma~\ref{r-part}. Hence $g>18$, which contradicts Table~2. Assume that $a\in\{3,4,6\}$ and suppose that $\alpha>\beta+1$. Similarly to
the previous case, $(q^c-1)_{\{r\}}=(q^b-1)_{\{r\}}\cdot{f}_{\{r\}}$, so $f_{\{r\}}\geq{r^2}$ and $c>18$; a contradiction.

Suppose that $r=3$. Then $a,b\in\{1,2\}$. If $a=1$ then
$3^{\alpha+1}~|~(u-1)(u^6+u^3+1)/(u-1,2)\in\omega(S)$, and if $a=2$ then $3^{\alpha+1}~|~(u+1)(u^6-u^3+1)/(u-1,2)\in\omega(S)$ (see Tables~\ref{tab:table2} and~\ref{tab:table3}).  In both cases $3^{\alpha+1}\in\omega(L)$. Therefore, if $g$ is the least positive integer such that $3^{\alpha+1}$ divides $q^g-1$, then $g=by\leq18$ for some positive integer~$y$. Since $b=1$ or $b=2$ and $(q^9+1,3)=1$ in the former case, an application of Lemma~\ref{r-part} gives $y=9$. So $L$ should contain a semisimple element whose $3$-part is equal to $(q^9-1)_{\{3\}}$ for $b=1$ and to $(q^9+1)_{\{3\}}$ for $b=2$. Looking over Table~2, we find that this is impossible.

Let now $r=2$. Then $a,b\in\{1,2\}$. Since $u^4-1\in\omega(S)$, Lemma \ref{r-part} implies that $(u^4-1)_{\{2\}}=2^{\alpha+2}$, so $2^{\alpha+2}\in\omega(L)$.
Let $g$ be a least positive integer such that $2^{\alpha+2}$ divides $q^g-1$.
If $g$ is odd, then $(q^g-1)_{\{2\}}=(q-1)_{\{2\}}\leq2^{\beta}<2^\alpha$, so $c$ is even. Choose $\tau\in\{+,-\}$ such that $q\equiv\tau1\pmod4$. Now Lemma~\ref{r-part} implies that $(q^g-1)_{\{2\}}=g_{\{2\}}\cdot(\tau{q}-1)_{\{2\}}\leq{g}_{\{2\}}\cdot2^{\alpha-1}$. Therefore, $g$ is divisible by $8$. Similarly to the case $r=3$, we conclude that $(q^8-1)_{\{2\}}\in\omega(L)$ and derive a contradiction using Table~2.
\end{prf}

\begin{lemma}\label{k_i_estimation}
Let $n$ be an integer and $n\geq2$. Then
\newline $(1)$ $k_1(n)k_2(n)=(n^2-1)/(2,n-1)$ and $n^2/4\leq{k}_1(n)k_2(n)\leq{n^2};$
\newline $(2)$ $k_3(n)k_6(n)=(n^4+n^2+1)/(3,n^2-1)$ and $n^4/3\leq{k}_3(n)k_6(n)\leq(5/4)n^4;$
\newline $(3)$ $k_4(n)=(n^2+1)/(2,n-1)$ and $n^2/2\leq{k}_4(n)\leq(5/4)n^2;$
\newline $(4)$ $k_5(n)k_{10}(n)=(n^8+n^6+n^4+n^2+1)/(5,n^2-1)$ and $n^8/5\leq{k}_5(n)k_{10}(n)\leq(4/3)n^8;$
\newline $(5)$ $k_7(n)=(n^6+n^5+n^4+n^3+n^2+n+1)/(7,n-1)$, $k_{14}(n)=(n^6-n^5+n^4-n^3+n^2-n+1)/(7,n+1)$, and $n^{12}/7\leq{k}_7(n)k_{14}(n)\leq(3/2)n^{12};$
\newline $(6)$ $k_8(n)=(n^4+1)/(n-1,2)$ and $(n^4)/2\leq{k}_8(n)\leq(17/16)n^4;$
\newline $(7)$ $k_9(n)=(n^6+n^3+1)/(3,n-1)$, $k_{18}(n)=(n^6-n^3+1)/(3,n+1)$, and $n^{12}/3\leq{k}_9(n)k_{18}(n)\leq(65/64)n^{12};$
\newline $(8)$ $k_{12}(n)=n^4-n^2+1$ and $(3/4)n^4\leq{k}_{12}(n)\leq{n^4};$
\newline $(9)$ $k_{15}(n)=n^8+n^7-n^5-n^4-n^3+n+1$, $k_{30}(n)=n^8-n^7+n^5-n^4+n^3-n+1$, and $(3/4)n^{16}\leq{k}_{15}(n)k_{30}(n)\leq{n^{16}};$
\newline $(10)$ $k_{20}(n)=(n^8-n^6+n^4-n^2+1)/(5,n^2+1)$ and $(4/25)n^8\leq{k}_{20}(n)\leq{n}^8;$
\newline $(11)$ $k_{24}(n)=n^8-n^4+1$ and $(15/16)n^8\leq{k}_{24}(n)\leq{n^8}.$
\end{lemma}

\begin{prf} The lemma is a direct consequence of the formula~(\ref{eq:ki}), Lemma~\ref{cyclomatic_polynomials} and straightforward computations.
\end{prf}

\begin{lemma} $S\not\simeq E_8(u)$.
\end{lemma}
\begin{prf}
Assume the contrary and let $S\simeq{E_8(u)}$. Since $t(v,S)=5$ (see Figure~2), it follows that $t(v,L)=t(v,G)\leq6$ by Lemma~\ref{l:Vasiliev theorem}, so  $v\in\{p\}\cup{R_1(q)}\cup{R_2(q)}$ (see Figure~1). Suppose that $q=5,7,9,11,13,$ or $17$. Then $v=2,3,5,7,$ or $p$. If $v=p$ then $r_{30}(p)\in\omega(G)\setminus\omega(L)$. Let $v=2$. In this case
$41\in{R}_{20}(2)$ and $31\in{R}_5(2)$ lie in $\pi(S)$, but $e(41,5)=20$, $e(41,7)=e(41,11)=e(41,13)=e(41,17)=40$, and $e(31,9)=15$, so
either $41\in\omega(G)\setminus\omega(L)$, or $31\in\omega(G)\setminus\omega(L)$; a contradiction.
If $v=3$ then $4561\in{R}_{15}(3)$, however $e(4561,5)=190$, $e(4561,7)=2280$, $e(4561,9)=15$, and  $e(4561,11)=e(4561,13)=e(4561,17)=4560$.
Therefore $4561\in\omega(G)\setminus\omega(L)$; a contradiction. Suppose $v=5$, then either $q=9$, or $q=11$. Note that $1741\in{R_{15}(5)}$ and $e(1741,9)=e(1741,11)=435$, so $1741\in\omega(G)\setminus\omega(L)$.
If $v=7$ then $q=13$. Since $31\in{R}_{15}(7)$ and $e(31,13)=30$, we get $31\in\omega(G)\setminus\omega(L)$; a contradiction.
Thus, we may assume that $q>17$.

Lemma~\ref{divisible} implies that $k_9(\varepsilon{q})$ divides $k_i(u)$ for some $i\in\{15,20,24,30\}$. It follows from Lemma~\ref{k_i_estimation} that $k_9(\varepsilon{q})\geq{(q^6-q^3+1)/3}\geq(99/300){q}^6$ and $k_i(u)\leq{u^8+u^7-u^5-u^4-u^3+u+1}\leq(4/3)u^8$. Therefore $q^6\leq(400/99)u^8$.

Set $I=\{1,2,3,4,5,6,7,8,9,10,12,14,15,18,20,24,30\}$ and $J=\{1,2,3,4,5,6,7,8,9,10,12,14,18\}$. Then $\bigcup\limits_{i\in{I}}R_i(u)\subseteq\pi(S)\subseteq\pi(L)=\{p\}\cup(\bigcup\limits_{j\in{J}}R_j(q))$. Put $a=\prod\limits_{i\in{I}}k_i(u)$ and $b=\prod\limits_{j\in{J}}k_j(q)$.
If $p\in{R_j(u)}$ for some $j$ and $p^\alpha$ divides $k_j(u)$, then $p^\alpha\in\omega(L)$. Since $q>17$ and $q$ is a $p$-power, Lemma \ref{unipotent} yields $p^\alpha\leq{q}$. Therefore, $a$ divides $35\cdot{b}\cdot{q}$ due to Lemma~\ref{R-inclusions}.
Lemma~\ref{k_i_estimation} implies that $b\leq{q^2}\cdot\frac{5q^4}{4}\cdot\frac{5q^2}{4}\cdot\frac{4q^8}{3}\cdot\frac{3q^{12}}{2}\cdot\frac{17q^4}{16}
\cdot\frac{65q^{12}}{64}\cdot{q^4}<
\frac{5\cdot5\cdot4\cdot3\cdot17\cdot65}{4\cdot4\cdot3\cdot2\cdot16\cdot64}q^{48}<(7/2)q^{48}$ and
$a\geq\frac{u^2}{4}\cdot\frac{u^4}{3}\cdot\frac{u^2}{2}\cdot\frac{u^8}{5}\cdot\frac{u^{12}}{7}\cdot\frac{u^4}{2}\cdot\frac{u^{12}}{3}
\cdot\frac{3u^4}{4}\cdot\frac{3u^{16}}{4}\cdot\frac{4u^8}{25}\cdot\frac{15u^8}{16}
\geq\frac{3\cdot3\cdot4\cdot15}{4\cdot3\cdot2\cdot5\cdot7\cdot2\cdot3\cdot4\cdot4\cdot25\cdot16}u^{80}>(1/59734)u^{80}$.
It follows that $(1/59734)u^{80}<35\cdot(7/2)q^{49}$, hence $u^{80}<7400000q^{49}$. On the other hand, $q^6\leq(400/99)u^8$, so
$q^{60}\leq(400/99)^{10}u^{80}<1400000u^{80}$. Therefore $q^{60}<(1400000\cdot7400000)q^{49}<17^5\cdot17^6q^{49}$ and so $q^{11}<17^{11}$, whence $q<17$; a contradiction.
\end{prf}

\begin{lemma} If $S\simeq E_7(u)$, then $u=q$.
\end{lemma}
\begin{prf}
Assume that $S\simeq E_7(u)$ and $u\neq{q}$.  If $r\in{R}_1(u)\cup{R}_2(u)$, then $t(r,S)=3$ and so $t(r,L)=t(r,G)\leq{t}(S)+1=4$.
Therefore, $r\in{R}_1(q)\cup{R}_2(q)$. Lemma \ref{R-inclusions} implies that $k_1(u)\cdot{k}_2(u)$ divides $k_1(q)\cdot{k}_2(q)=(q^2-1)/2$.

Suppose that $v=2$. Then $k_1(u)\cdot{k}_2(u)=u^2-1$ is odd. Since $8~|~(q^2-1)$, we obtain $q^2-1\geq8(u^2-1)\geq(2u)^2$, so $q>2u$.

Let now $v\neq{2}$. Choose $\tau\in\{+,-\}$ such that $2$ is not adjacent with
$r_7(\tau{u})$ and $r_9(\tau{u})$ in $GK(S)$ (see Table \ref{tab:table1}). Then $R_7(\varepsilon{q})\subseteq{R}_i(\tau{u})$ and $R_7(\varepsilon{q})\subseteq{R}_j(\tau{u})$, where $i$ and $j$ are distinct numbers from $\{7,9\}$, due to Lemmas~\ref{l:Vasiliev theorem} and~\ref{2coclique}, and Table \ref{tab:table1}.
Lemma~\ref{divisible} implies that one of the numbers $(u^7-\tau1)/2$ and $(u-\tau1)(u^6+\tau{u}+1)/2$ divides $(q^7-\varepsilon1)/2$, while
the other divides $(q-\varepsilon1)(q^6+\varepsilon{q}+1)/2$. In particular, the greatest common divisor of these numbers divides the greatest common divisors of $(q^7-\varepsilon1)/2$ and $(q-\varepsilon1)(q^6+\varepsilon{q}+1)/2$. Therefore, $u-\tau1$ divides $q-\varepsilon1$ and so $q-\varepsilon1=l(u-\tau1)$ for some positive integer~$l$. On the other hand, $k_1(u)\cdot{k}_2(u)=(u^2-1)/2$ divides $(q^2-1)/2$ and we take a positive integer $k$ such that $q^2-1=k(u^2-1)$. Since $q\neq u$, we have $q>u$, so $k>1$.

Suppose that $k<4$. Then either $k=2$ or $k=3$. If $l>k$, then $q+1\geq{q}-\varepsilon1>2(u-\tau1)\geq2u-2$, so $(q+2)/2\geq{u}$.
On the other hand, $q^2-1=k(u^2-1)<(q-\varepsilon1)(u+\tau1)$. It follows that $u+\tau1>q+\varepsilon1$, hence $u+1>q-1$. This implies that $(q+2)/2\geq{u}\geq{q-1}$, which is impossible due to $q>4$. Observe that $l=k$ implies $q+\varepsilon1=u+\tau1$, so the cases
$l=k$ and $l=1$ are the same (it is sufficient to replace $\tau$ on $-\tau$ and $\varepsilon$ on $-\varepsilon$). Assume that $l=1$ and $u-\tau1=q-\varepsilon1$. Since $q>u$, it follows that $q-1=u+1$. Then $k(u-1)=q+1$. Therefore, $2=q+1-(q-1)=k(u-1)-(u+1)=(k-1)u-(k+1)$. Hence $u=(k+3)/(k-1)=1+4/(k-1)$. It follows that $(u,q)\in\{(3,5),(5,7)\}$.
If $(u,q)=(5,7)$, then $5^3-1=4\cdot31$ and $e(31,7)=15$, so $31\in\pi(S)\setminus\pi(L)$; a contradiction. If $(u,q)=(3,5)$, then $3^4+1=2\cdot41$ and $e(41,5)=20$, and we derive a contradiction because $41\in\pi(S)\setminus\pi(L)$. Thus, $1<l<k<4$, hence $l=2$ and $k=3$.
This yields $2(q+\varepsilon1)=3(u+\tau1)$, so $3(q-\varepsilon1)=6(u-\tau1)$ and $4(q+\varepsilon1)=6(u+\tau1)$.
Therefore, $q+\varepsilon7=\tau12$. Hence either $q=5$ or $q=19$.
If $q=5$ then $u=3$ which is impossible as proved above. If $q=19$, then $u=11$ and $61\in\pi(S)\setminus\pi(L)$; a contradiction.
Thus we may assume that $k\geq4$ and $q^2-1\geq4(u^2-1)$. Straightforward calculations show that $q>3u/2$ in this case. Thus, we always have the inequality $q>3u/2$.

The inequality $q>3u/2$ implies that $3\cdot{k}_9(\varepsilon{q})\geq{q^6}+\varepsilon{q^3}=q^3(q^3+\varepsilon1)>(3{u}/2)^3((3u/2)^3+\varepsilon1)>11u^6-4u^3>3\cdot(u^6+u^5+u^4+u^3+u^2+u+1)\geq3\cdot{k}_i(u)$, where $i\in\{7,9,14,18\}$ due to Lemma~\ref{k_i_estimation}. On the other hand, Lemma \ref{divisible} implies that $k_9(\varepsilon{q})$ divides one of $k_i(u)$, where $i\in\{7,9,14,18\}$; a contradiction. Thus $u=q$. The lemma is proved, which completes the proof of the theorem.
\end{prf}

\end{document}